# О КВАЗИПЕРИОДИЧЕСКИХ РАЗБИЕНИЯХ, ОБЛАДАЮЩИХ СИММЕТРИЕЙ ПЯТОГО ПОРЯДКА

**Александр С. Прохода, к.ф.-м.н.**

Украина

Применен теоретико-групповой подход к построению квазипериодических разбиений евклидовой плоскости, обладающих симметрией пятого порядка. Из бесконечного множества вариантов квазипериодических разбиений плоскости, обладающих диэдральной группой симметрии $D_5$, особое внимание уделено тем разбиениям, которые можно получить одним универсальным набором, в состав которого входит пять различных плиток. Определены геометрические характеристики тайлов из этого набора. Показано, что данным набором плиток можно осуществлять топологически различные разбиения плоскости, обладающие поворотными симметриями, как пятого, так и десятого порядка.

___________________________

e-mail: a-prokhoda@mail.ru

# I. ВВЕДЕНИЕ

Иоганн Кеплер в 1619 г. в своем трактате «HARMONICES MUNDI» комплексно исследовал мозаики на плоскости и начал толковать их с математической точки зрения. В частности Кеплер изучил и построил мозаики из выпуклых и не выпуклых многоугольников (звезд). Роджер Пенроуз в [1] установил геометрические характеристики двух ромбов и правила их сочетания, при помощи которых можно выполнить квазипериодическое разбиение (P3) евклидовой плоскости. Рисунок, приведенный Кеплером в трактате «Гармония мира» на стр. 61 (если не обращать внимания на тени, которые очевидно символизируют об объемности фигуры) совпадает с центральной частью мозаики Пенроуза (P3). Кеплер в первой главе, также акцентирует внимание на невозможность выполнить разбиение плоскости грань в грань конгруэнтными пентагонами (гептагонами). Однако, несколькими многоугольниками одновременно, в том числе и звездчатыми, гениальному Кеплеру удалось воспроизвести фрагменты квазипериодических разбиений обладающих поворотными симметриями 5 и 8 порядков (см. рисунки в конце первой главы).

Шехтман в [2] экспериментально установил квазикристаллическую структуру в сплаве $Al_{86}Mn_{14}$ (данный квазикристаллический сплав впоследствии получил название «шехтманит»). Электроннограммы приведенные Шехтманом в [2] свидетельствовали о наличии икосаэдральной симметрии в отдельных зернах данного материала.

В одной из математических моделей квазикристаллы рассматриваются как проекции многомерных кристаллов. Дышлис в [3, 4] подробно описал этот подход.

Рассмотрим целочисленную решетку, построенную на 5 векторах единичной длины, попарно ортогональных один к одному в пятимерном евклидовом пространстве $E^5$. Эти вектора образуют ортонормированный базис пространства $E^5$.

Очевидно, что циклическая подстановка этих векторов по правилу $e_1 \to e_2 \to e_3 \to e_4 \to e_5 \to e_1$ есть операцией симметрии пятого порядка. Спроектировав эти вектора на пространство $E^3$ или $E^2$ получим совокупность векторов которые также имеют симметрию пятого порядка. Также можно рассматривать одномерное пространство $E^1$ натянутое на сумму векторов $e_1 + e_2 + e_3 + e_4 + e_5$ и ортогональное дополнение к нему четырехмерного евклидова пространства $E^5 = E^4 \oplus E^1$. Проектируя четырехмерную решетку пространства $E^4$ с симметрией пятого порядка, можно получить квазирешетку состоящую из 4-х векторов и имеющую, также симметрию пятого порядка (поскольку $E^1$ и $E^5$ не изменяются при перестановке векторов $e_1, e_2, e_3, e_4, e_5$).

Приведем понятие $n$-мерной решетки.

**Определение 1.** Система векторов $\varepsilon_1, \varepsilon_2, \ldots, \varepsilon_k$ евклидовой плоскости называется целочисленно-линейно независимой, если ее целочисленная линейная комбинация $z_1\varepsilon_1 + z_2\varepsilon_2 + \ldots + z_k\varepsilon_k = 0$ возможна только тогда, когда все коэффициенты $z_i = 0$ при $i = 1, 2, \ldots, k$. Если же $z_1\varepsilon_1 + z_2\varepsilon_2 + \ldots + z_k\varepsilon_k = 0$ при некотором $z_i \neq 0$, то система называется целочисленно-линейно зависимой.

В каждой целочисленно-линейно зависимой системе векторов существует целочисленно-линейно независимая подсистема векторов.

Приведем пример целочисленно-линейно независимой системы векторов, линейная комбинация которых имеет симметрию пятого порядка. Для построения такой системы рассмотрим все корни пятой степени из 1. Они все удовлетворяют уравнению $x^5 = 1$ или $x^5 - 1 = 0$, которое раскладывается на множители следующим образом $(x - 1)(x^4 + x^3 + x^2 + x + 1) = 0$. Число $1$ является корнем уравнения $x - 1 = 0$, а все другие корни будут корнями уравнения $x^4 + x^3 + x^2 + x + 1 = 0$, которое неприводимо над полем рациональных чисел $Q$, потому что все корни кроме $1$ комплексные и являются степенями $\cos(72°) + i\sin(72°)$. Вспомним, что $\cos(72°) = (\sqrt{5} - 1)/4$ и $2\cos(36°) = 2\sin(54°) = (1 + \sqrt{5})/2 = \tau$. Система всех $5$ векторов $1, \varepsilon, \varepsilon^2, \varepsilon^3, \varepsilon^4$ целочисленно-линейно зависима (по теореме Виета $1 + \varepsilon + \varepsilon^2 + \varepsilon^3 + \varepsilon^4$ равняется коэффициенту при $x^4$ в уравнении $x^5 - 1 = 0$, то есть равняется $0$). Если же мы возьмем, числа $1, \varepsilon, \varepsilon^2, \varepsilon^3$ то эта система целочисленно-линейно независима. Действительно допустим, что она целочисленно-линейно зависима, то есть удовлетворяет уравнению с целыми коэффициентами степени $3$.

$$z_0 + z_1\varepsilon + z_2\varepsilon^2 + z_3\varepsilon^3 = 0.$$

Этого не может быть, потому что $\varepsilon$ удовлетворяет неприводимому над полем $Q$ уравнению четвертой степени, а мы получили, что $\varepsilon$ удовлетворяет уравнению третьей степени с целыми коэффициентами.

Множество $z_0 + z_1\varepsilon + z_2\varepsilon^2 + z_3\varepsilon^3$, где $z_i$ пробегает все целые числа называется квазирешеткой ранга $4$. Покажем, что она совмещается с собой при поворотах на углы кратные $72°$. Действительно при умножении комплексных чисел их аргументы складываются. Поэтому оператор поворота является оператором умножения на число

$$\varepsilon = \exp(2\pi i/5) = \cos(2\pi/5) + i\sin(2\pi/5).$$

Если умножить $z_0 + z_1\varepsilon + z_2\varepsilon^2 + z_3\varepsilon^3$ при фиксированных $z_0, z_1, z_2, z_3$ на $\varepsilon$ имеем

$$z_0\varepsilon + z_1\varepsilon^2 + z_2\varepsilon^3 + z_3\varepsilon^4. \qquad (1)$$

Однако система $1, \varepsilon, \varepsilon^2, \varepsilon^3, \varepsilon^4$ целочисленно-линейно зависима, то есть один из векторов является целочисленной линейной комбинацией остальных

$$\varepsilon^4 = -1 - \varepsilon - \varepsilon^2 - \varepsilon^3.$$

Подставив $\varepsilon^4$ в (1) имеем

$$z_0\varepsilon + z_1\varepsilon^2 + z_2\varepsilon^3 + z_3(-1 - \varepsilon - \varepsilon^2 - \varepsilon^3) =$$
$$= -z_3 + (z_0 - z_3)\varepsilon + (z_1 - z_3)\varepsilon^2 + (z_2 - z_3)\varepsilon^3 =$$
$$= z_0'1 + z_1'\varepsilon + z_2'\varepsilon^2 + z_3'\varepsilon^3,$$

где $z_0' = -z_3$, $z_1' = z_0 - z_3$, $z_2' = z_1 - z_3$, $z_3' = z_2 - z_3$. То есть квазирешетка совмещается с собой при повороте на угол $72°$.

Система корней дает возможность найти группу симметрии системы корней и наоборот, если задана группа Вейля, то действуя этой группой на базисные векторы можно получить систему корней.

**Определение 2.** Если $\Sigma_n$ – система корней векторного пространства $V$, то группой Вейля системы $\Sigma_v$ называется подгруппа $W$ в полной линейной группе $GL(V)$, порожденная отражениями $S_{\vec{a}}$, $\vec{a} \in \Sigma_n$.

Группа $W$ является, очевидно, подгруппой группы $\mathrm{Aut}(\Sigma_n)$ автоморфизмов пространства $V$, которая совмещает множество $\Sigma_n$ с собой.

Приведем более общее определение системы корней, которое дает возможность классифицировать группы симметрий правильных многогранников.

**Определение 3.** Конечный набор векторов $\Phi \subset \mathrm{E}^n$ называется системой корней если он удовлетворяет следующим двум аксиомам:

R1: $\Phi \cap \mathrm{R}_a = \{a, -a\}$;

R2: $S_a\Phi = \Phi, \forall\, a \in \Phi$.

Здесь $\{a, -a\}$ фиксированная пара векторов (единичной длины).

По аксиоматически определенной системе корней $\Phi$ можно построить группу отражений $W$: это будет группа порожденная отражениями всех корней из $W$.

**Определение 4.** Пусть $O(n, R)$ ортогональна группа пространства $E^n$. Конечная подгруппа $W \in E^n$ называется группой порожденной отражениями, если существуют такие отражения $S_{a_1}, \ldots, S_{a_n}$, которые являются образующими элементами этой группы.

Бесконечные группы также могут быть порождены отражениями, потому что группа движений пространства $E^n$ порождена, как известно отражениями.

Из определения принципа симметрии элементов симметрии следует следующая теорема.

**Теорема 1.** Пусть $t \in O(n, R)$ произвольное ортогональное преобразование является $S_a$ отражением относительно вектора $a$ (который является корнем). Тогда $tS_a t^{-1} = S_{t_a}$ – отражение относительно вектора $t_a$ в частности, если $w \in W$ то $S_{w_a} \in W$.

↑ Очевидно, что $tS_a t^{-1}(ta) = tS_a(a) = t(-a) = -ta$, тогда $tS_a t^{-1}$ переводит вектор $ta$, в противоположный ему вектор $-ta$.

Покажем, что $tS_a t^{-1}$ оставляет на месте плоскость $H_{ta}$ – зеркало отражения. Действительно условие, что вектор $\vec{b}$ принадлежит $H_a$ равносильно тому, что $tb$ принадлежит $H_{ta}$, так как $(a, b) = (ta, tb)$, поэтому $(tS_a t^{-1})(tb) = tS_a(tb) = tb$. ↓

Система корней ранга 2 позволяет построить решетку основными векторами, которой будут простые корни $\vec{a}$ и $\vec{b}$ этой системы. Она называется корневой решеткой ранга 2. Как абстрактная группа – это свободная абелева группа ранга 2. Особенностью корневой решетки является то, что ей принадлежат все корни системы корней.

**Определение 5.** Если $\Sigma_n$ – система корней ранга 2 пространства $E^2$, то дискретная группа, порожденная отражениями от прямых

$L_{a,k} = \{x \in E^2 \mid \frac{2(\vec{a}\vec{x})}{(\vec{a}\vec{a})}\vec{a}\} = k, k \in Z$ называется аффинной группой Вейля системы корней.

Поскольку решетка связана с группой параллельных переносов (подгруппа $T$ параллельных переносов, которая сохраняет решетку), то группа Вейля $W_{\vec{a}}(\Sigma_2)$ раскладывается в полупрямое произведение группы Вейля $W(\Sigma_2)$ системы корней на $T$, то есть

$$W_{\vec{a}}(\Sigma_2) = T \lambda W(\Sigma_2).$$

Напомним, что полная группа симметрии решетки называется группой Браве. Отсюда следует что аффинная группа Вейля системы корней является либо группой Браве либо ее подгруппой. Группа Браве корневой решетки построенная на простых корнях $\vec{a}$ и $\vec{b}$ является аффинной группой Вейля.

В топологии и в современной геометрии широко распространены специфические операции, которые позволяют из одних многообразий получать другие многообразия. Это операции склеивания многообразий и обратная к ней операция разрезания многообразий, операция приклеивания и переклеивания, операция заплаты и образования дыр. Все эти операции получили название хирургии многообразий [5, 6]. Великим "Хирургом" был Уильям Тёрстон, который разработал метод исследования трехмерных многообразий, основанный на разрезании их на куски, допускающий локально-однородную метрику. В настоящей работе развивается также идея, принадлежащая Тёрстону высказанная в [7], получения модели идеального кристалла евклидовой геометрии. Строится алгоритм, позволяющий по схеме склейки фундаментальной области фундаментальной группы многообразия, получать разбиения носителя геометрии на ячейки, декорируемые атомами (модель идеального кристалла). В данном исследовании используется тот факт, что модель идеального или реального кристалла можно получить путем действия квазикристаллографической группы на фундаментальную область этой группы, причем группа задается с помощью ее генетического кода. С другой стороны схема склейки кодируется, словом, принадлежащим фундаментальной группе многообразия, с помощью которого получается разбиение накрытия многообразия.

## II. АЛГОРИТМ ПОСТРОЕНИЯ КВАЗИРЕШЕТОК И СООТВЕТСТВУЮЩИХ РАЗБИЕНИЙ ЕВКЛИДОВОЙ ПЛОСКОСТИ

Проанализируем квазирешетку ранга 4, как свободную абелевую группу ранга 4 построенную на корнях пятой степени из единицы с группой Вейля $D_5$, порядок которой равен 10. Отметим, что характерные плоские сечения материалов с икосаэдральной группой симметрии, имеют симметрию пятого порядка и в качестве их моделей служат в том числе и разбиения Пенроуза де Брюина плоскости $E^2$.

Рассмотрим множество $T \otimes Q \subset R^4$ рациональных линейных комбинаций векторов, которые принадлежат $T$, которое является двумерным подпространством пространства $R^2$ над полем $Q_{[\sqrt{5}]}$. Числа поля $Q_{[\sqrt{5}]}$ можно подать в виде $p + q\sqrt{5}$, где $p$, $q$ – рациональные числа, а целые алгебраические числа этого поля имеют вид $m + n\tau$, где $m$, $n$ – целые числа, при этом $n$ должно быть парным числом: $n = 2n_1$. Таким образом, имеем $m + n\tau = m + n_1(1 + \sqrt{5}) = m + n_1 + n_1\sqrt{5}$, причем это число и число $m + n_1 - n_1\sqrt{5}$ должны удовлетворять уравнению

$$x^2 - 2(m + n_1)x - (m + n_1)^2 - 5n_1^2 = 0 \text{ или}$$
$$x^2 - (2m + n)x - (m^2 + mn - n^2) = 0.$$

Приходим к выводу, что четырехмерное действительное векторное пространство можно рассматривать как двумерное подпространство над полем $Q_{[\sqrt{5}]}$ для некоторого натурального $d$ (в данном случае $d = 5$).

Поскольку квазирешетка ранга 4 может рассматриваться, как корневая решетка обобщенной системы корней ранга 2, которая построена на простых корнях, то далее найдем обобщенную систему корней ранга 2 с группой Вейля $D_5$.

Проведем зеркала отражений – оси симметрии пентагона, они проходят через его вершины и середины противоположных сторон и делят их пополам. Таким образом, имеем 10 корней – пять положительных и пять отрицательных, которые являются корнями десятой степени из 1. Приведем их к общему началу, воспользовавшись тем свойством, что вектора свободные. Один из них, это вектор $\vec{a}=1$, угол между соседними векторами равен 36°. На роль простых корней претендуют: корень $\exp(2\pi i/10) = \cos(144°) + i\sin(144°)$, угол между векторами $\vec{a}$ и $\vec{b}$ равен 144° он равен первообразному корню пятой степени из 1 $\varepsilon = \cos(72°) + i\sin(72°)$ в квадрате $\varepsilon^2$ или вектору $\vec{b} = \varepsilon_1^4$, где $\varepsilon_1 = \cos(36°) + i\sin(36°)$, $\varepsilon = \varepsilon_1^2$. Поскольку $\cos(36°) = (1+\sqrt{5})/4 = \tau/2$, то выберем вектор $\vec{a}=1$, а вектор $\vec{b} = \exp(\pi i/5) = \cos(144°) + i\sin(144°)$. Чтобы определить остальные корни воспользуемся формулой для отражения $\hat{S}_{\vec{a}}(\vec{b}) = \vec{b} - \frac{2(\vec{a},\vec{b})}{(\vec{a},\vec{a})}\vec{a} = \vec{b} + 2\cos(36°)\vec{a}$, $\hat{S}_{\vec{b}}(\vec{a}) = \vec{a} - \frac{2(\vec{b},\vec{a})}{(\vec{b},\vec{b})}\vec{b} = \vec{a} + \tau\vec{b}$. Чтобы найти $\hat{S}_{\vec{a}}\hat{S}_{\vec{b}}(\vec{a})$ воспользуемся тем, что $\tau^2 = \tau + 1$ или $x^2 - x - 1 = 0$ уравнением корнем, которого является $\tau$.

Далее

$(\hat{S}_{\vec{a}}\hat{S}_{\vec{b}})\vec{a} = \hat{S}_{\vec{a}}(\vec{a}+\tau\vec{b}) = -\vec{a} + \tau\hat{S}_{\vec{a}}(\vec{b}) = -\vec{a} + \tau(\vec{b}+\tau\vec{a}) = -\vec{a} + \tau\vec{b} + \tau^2\vec{a} =$

$= -\vec{a} + \tau\vec{b} + \vec{a} + \tau\vec{a} = \tau(\vec{a}+\vec{b})$. Таким образом, имеем пять положительных корней: $\vec{a}, \vec{b}, \vec{a}+\tau\vec{b}, \vec{b}+\tau a, \tau(\vec{a}+\vec{b})$. К ним необходимо прибавить столько же отрицательных корней. Окончательно имеем

$$\Sigma_5 = \{\pm\vec{a}, \pm\vec{b}, \pm(\vec{a}+\tau\vec{b}), \pm(\vec{b}+\tau a), \pm\tau(\vec{a}+\vec{b})\}.$$

Далее выпишем матрицы отражений в координатах базиса $\vec{a}, \vec{b}$. Имеем $\hat{S}_{\vec{a}}(\vec{a}) = -\vec{a} = -1\vec{a} + 0\vec{b}$, $\hat{S}_{\vec{a}}(\vec{b}) = \vec{b} + \tau\vec{a}$ поэтому $\hat{S}_{\vec{a}}$ имеет матрицу $S_a = \begin{pmatrix} -1 & 0 \\ \tau & 1 \end{pmatrix}$ определитель которой равен, как и должно быть –1, $S_a^2 = \begin{pmatrix} 1 & 0 \\ 0 & 1 \end{pmatrix}$,

$\hat{S}_{\vec{b}}(\vec{a}) = \vec{a} + \tau\vec{b}$, $\hat{S}_{\vec{b}}(\vec{b}) = -\vec{b}$, $S_b = \begin{pmatrix} 1 & \tau \\ 0 & -1 \end{pmatrix}$, $S_b^2 = \begin{pmatrix} 1 & 0 \\ 0 & 1 \end{pmatrix}$. Перемножив матрицы $S_aS_b$ получим матрицу поворота $S_aS_b = \begin{pmatrix} -1 & 0 \\ \tau & 1 \end{pmatrix}\begin{pmatrix} 1 & \tau \\ 0 & -1 \end{pmatrix} = \begin{pmatrix} -1 & -\tau \\ \tau & \tau \end{pmatrix}$, поскольку $\tau^2 = \tau + 1$, то $\tau^2 - 1 = \tau$. Определитель ее равен $\tau^2 - \tau = 1$, как и должно быть. Далее $(S_aS_b)^2 = \begin{pmatrix} 1-\tau^2 & \tau-\tau^2 \\ -\tau+\tau^2 & 0 \end{pmatrix} = \begin{pmatrix} -\tau & -1 \\ 1 & 0 \end{pmatrix}$,

$(S_aS_b)^4 = \begin{pmatrix} \tau & \tau \\ -\tau & -1 \end{pmatrix}$, $(S_aS_b)^5 = \begin{pmatrix} 1 & 0 \\ 0 & 1 \end{pmatrix}$. Таким образом, диэдральная группа $D_5$ задается $\hat{S}_a, \hat{S}_b$ двумя образующими и определяющими соотношениями $\hat{S}_a^2 = \hat{S}_b^2 = (\hat{S}_a\hat{S}_b)^5 = 1$. Кроме того, мы получили 7 матриц из 10. Найдем еще такие $S_aS_bS_a = \begin{pmatrix} -\tau & -\tau \\ 1 & \tau \end{pmatrix}$, $(S_aS_b)^3 = \begin{pmatrix} 0 & 1 \\ -1 & 0 \end{pmatrix}$, $S_bS_aS_b = \begin{pmatrix} \tau & 1 \\ -\tau & -\tau \end{pmatrix}$,

$(S_aS_b)^4S_a = \begin{pmatrix} 1 & \tau \\ 0 & -1 \end{pmatrix}$. Матрица $\begin{pmatrix} 1 & \tau \\ 0 & -1 \end{pmatrix}$ с определителем –1 является матрицей инверсии. Таким образом, $D_5$ можно определить при помощи образующих элементов и определяющих соотношений как: $\hat{S}_a \hat{S}_b = \hat{S}_c$,

$$i = \begin{pmatrix} 1 & \tau \\ 0 & -1 \end{pmatrix} = S_i, \qquad S_i S_c S_i = S_c^4, \qquad D_5 = \{S_c, S_i \mid S_c^5 = S_i^2 = 1\}, \qquad D_5 = (S_c)_5 \lambda S_i,$$

$$S_i S_c S_i = \begin{pmatrix} 1 & \tau \\ 0 & -1 \end{pmatrix} \begin{pmatrix} -1 & -\tau \\ \tau & \tau \end{pmatrix} \begin{pmatrix} 1 & \tau \\ 0 & -1 \end{pmatrix} = \begin{pmatrix} \tau & \tau \\ -\tau & -1 \end{pmatrix} = (S_a S_b)^4.$$ Видно, что все корни $\Sigma_5$ есть целыми алгебраическими комбинациями простых корней. Точечная группа $D_5$ является группой Браве квазикристалла.

## III. АНАЛИЗ НЕПЕРИОДИЧЕСКИХ РАЗБИЕНИЙ ЕВКЛИДОВОЙ ПЛОСКОСТИ, ОБЛАДАЮЩИХ СИММЕТРИЕЙ ПЯТОГО ПОРЯДКА

Известно [8], что если разрезать вдоль меньшей диагонали тонкий, а вдоль большей диагонали толстый ромбы Пенроуза, то в результате получатся равносторонние треугольники Робинсона с углами (36°, 72°, 72°) и (108°, 36°, 36°) соответственно. Далее треугольники Робинсона можно разрезать на меньшие, так, чтобы каждый из новых (меньших) был гомотетичен одному из исходных. Линейные размеры новых треугольников будут в $\tau$ раз меньше чем размеры исходных. Такое разрезание называется дефляцией. Обратное же преобразование (склеивание) называется инфляцией.

Так как фундаментальная область группы выбирается неоднозначно, то существует вариация в выборе тайлов в конкретном разбиении. На рис. 1 приведен фрагмент квазипериодического разбиения евклидовой плоскости переходящий в себя при поворотах на углы кратные 72° вокруг особой точки. Особая точка – геометрический центр множества. Очевидно, что данное разбиение обладает симметрией пятого порядка.

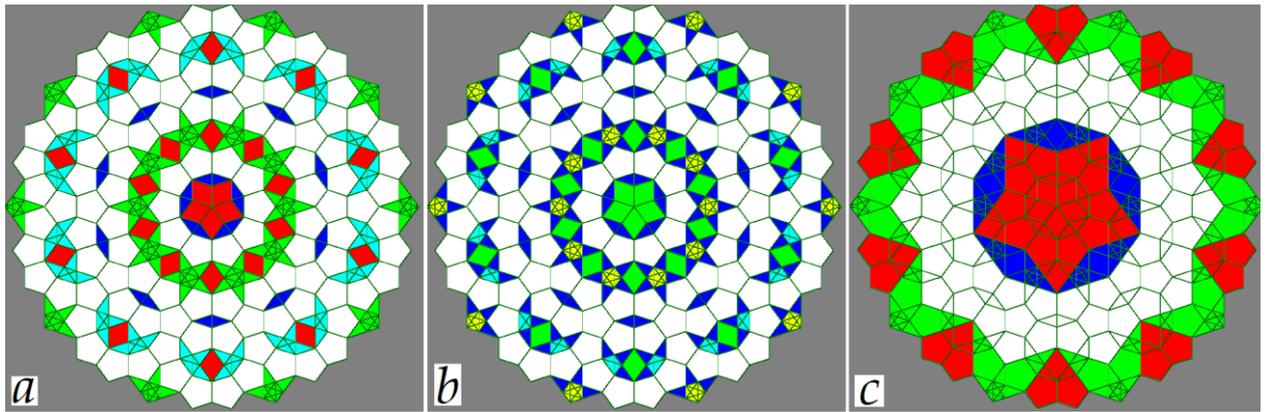

Рис. 1. Фрагмент квазипериодического разбиения евклидовой плоскости, обладающий симметрией пятого порядка:

на *a* и *b* тайлы разбиения выбраны по-разному;

*c* – визуализация операции гомотетического вращения разбиения *a*.

Фрагмент разбиения плоскости рис. 1, *a* выполнен следующими пятью геометрическими фигурами: толстым и тонким ромбами Пенроуза; пентагоном; пятиконечной звездой (пентаграммой); лодочкой Пенроуза. Пентаграмма образована из пентагона окруженного пятью треугольниками Робинсона (с углами 36°, 72°, 72°), основания которых приклеены к сторонам данного пентагона. Длина стороны пентагона, из которого получена пятиконечная звезда в τ раз меньше соответствующей длины пентагона указанного в перечне тайлов к разбиению рис. 1. Лодочка Пенроуза образована равнобедренной трапецией (с углами 72° и 108°) с приклеенными к ее сторонам и к меньшему основанию треугольниками Робинсона (с углами 36°, 72°, 72°). У данной трапеции меньшее основание и боковые стороны численно равны. Большее же основание данной трапеции в τ раз больше меньшего ее основания, т.е. равняется τ. Целое алгебраическое число τ целесообразно принять за единицу длины в разбиениях, исследованных в данной работе и вообще в любом разбиении обладающем симметрией пятого порядка. Также числу τ равняется длина стороны пентагона, длина стороны ромба Пенроуза и длина боковой стороной треугольника Робинсона.

На рис. 1, *b* приведено тоже разбиение, что и на рис. 1, *a,* однако пять различных тайлов разбиения выбраны несколько по иному принципу. Среди данных тайлов: треугольник Робинсона (с углами 36°, 72°, 72°), ромб Пенроуза (с углами 72° и 108°), равнобедренная трапеция (с углами 72° и 108°) и два подобных пентагона (коэффициент гомотетии между которыми равен τ). Заметим, что тайлы квазипериодического разбиения приведенного на рис. 1, *a* и рис. 1, *b* хоть и выбраны по-разному, однако их число в обоих случаях равно пяти. Отметим также, что возможны и другие варианты выбора тайлов. На рис. 1, *c* визуализирована операция гомотетического вращения осуществленная на разбиении рис. 1, *a*: коэффициент гомотетии равен $\tau^2 = \tau + 1$ и поворот на угол 72° вокруг особой точки. Для последующего этапа гомотетии линейные размеры плиток будут в $\tau^2$ раз больше в сравнении с соответствующими размерами закрашенных плиток изображенных на разбиении рис. 1, *c* т.е. в $\tau^4 = 3\tau + 2$ раз больше в сравнении с размерами плиток в разбиении рис. 1, *a*. Из рис. 1, *a* и рис. 1, *c* следует, что приведенные разбиения изоморфны.

Заметим, что в построенных разбиениях отношение количества конкретных типов тайлов сопоставимы с численными значениями некоторых степеней числа τ. Например, на разбиении рис. 1, *a* рассмотрим отношение числа пентагонов к числу толстых ромбов Пенроуза. Данное отношение стремится к $\tau^3$. Отношение числа пентагонов к числу тонких ромбов Пенроуза стремится к $\tau^4$; отношение числа пентагонов к числу звезд → $\tau^5$ и т.д.

На рис. 2, *a* приведено квазипериодическое разбиение плоскости, в котором участвует тот же набор плиток, что и в разбиении на рис. 1, *b*. Однако правила сочетания данных тайлов отличаются. Например, на рис. 1, *b* в особой точке сходятся пять толстых ромбов Пенроуза, а на рис. 2, *a* в особой точке сходятся десять треугольников Робинсона. На рис. 2 вершины, в которых сходятся тайлы, декорированы «атомами».

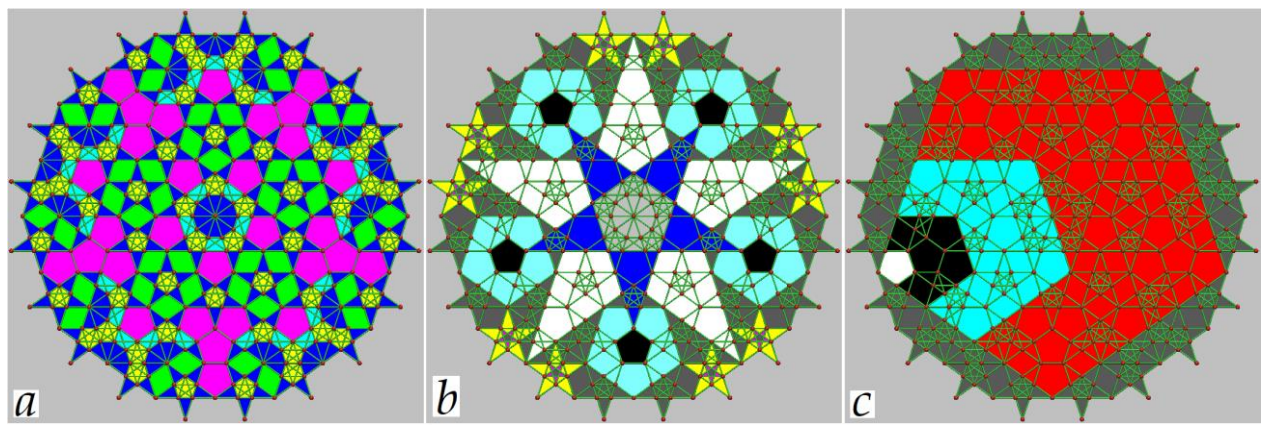

Рис. 2. Фрагмент квазипериодического разбиения плоскости *a*; рисунки *b* и *c* раскрашены таким образом, чтобы визуализировать в разбиении наличие различных фигур с симметрией пятого порядка.

На рис. 2, *b* визуализировано гомотетическое вращение относительно особой точки, меняющее пентаграмму на подобную ей пентаграмму. Геометрические центры пентаграмм совпадают, и совпадают с особой точкой разбиения. На рис. 2, *b* выделены пять пентагонов внутри каждого, из которых расположен пентагон меньшего размера (коэффициент подобия между которыми равен $\tau^2$), а также подобные звезды.

На рис. 2, *c* выделены гомотетические пентагоны, геометрический центр пентагона максимального размера совпадает с особой точкой разбиения. Не составляет большого труда определить коэффициенты гомотетии пентагона минимального размера с бо́льшими пентагонами.

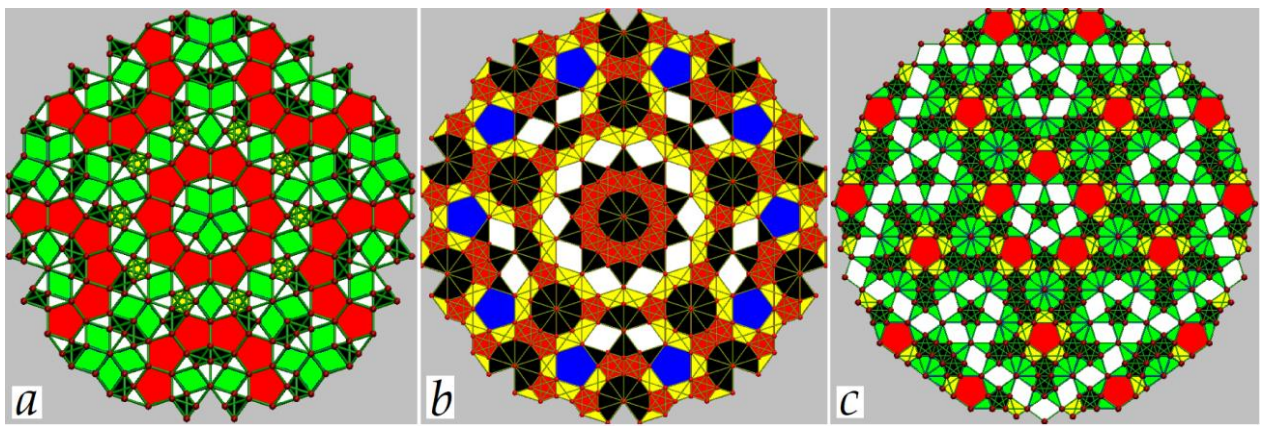

Рис. 3. Варианты квазипериодических разбиений плоскости, обладающие поворотной симметрией пятого порядка, выполненные одним и тем же набором тайлов.

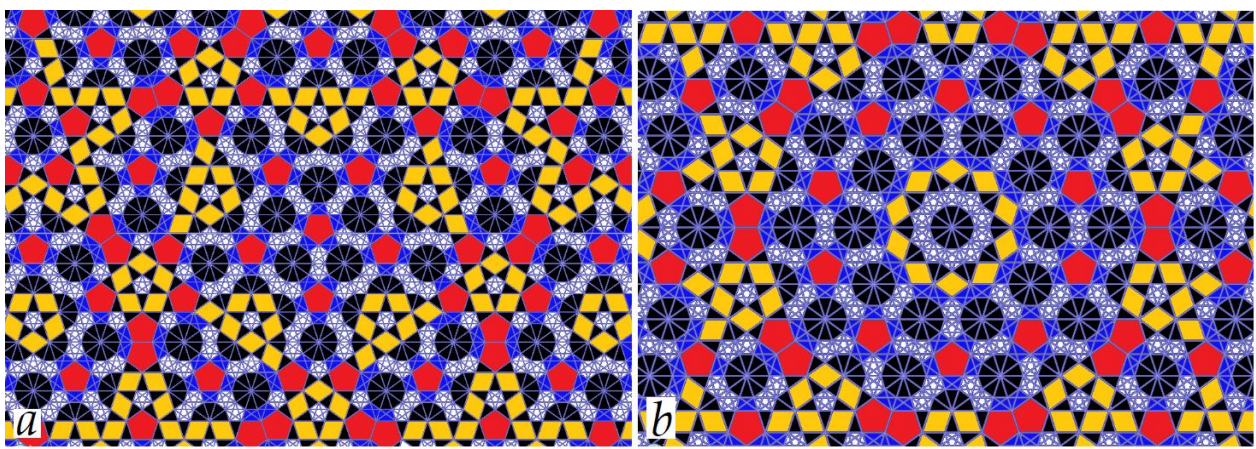

Рис. 4. Квазипериодические разбиения евклидовой плоскости; разбиение на *a* обладает поворотной симметрия пятого порядка, на *b* – десятого порядка.

В работе [9], мною уже были показаны варианты разбиения плоскости пятью тайлами (перечень тайлов тот же, что и в разбиениях, приведенных на рис. 1, *b* и рис. 2, *a*). Важно отметить что, таким набором плиток, возможно, произвести разбиение плоскости, совмещающееся с собой при поворотах на углы кратные 36° (т.е. разбиение плоскости, обладающее симметрией десятого порядка). На рис. 4 можно визуализировать этот факт.

На рис. 1, *b*, рис. 2, *a* и на рис. 3 приведены фрагменты квазипериодических разбиений плоскости, обладающие симметрией пятого порядка. Отметим, что данные разбиения хоть и получены при помощи одного набора плиток, но отношение, например числа пентагонов к числу трапеций в них разное. Однако эти отношения во всех случаях стремятся к численному значению различных, но определенных степеней числа τ. Центральная область рис. 3, *a* аналогична центральной области в разбиении Пенроуза (P3), т.е. 5 толстых ромбов и 5 тонких ромбов разбивают декагон. Однако декагон на рис. 3, *a* окружен с наружи десятью пентагонами. Отметим также сходство вблизи центральных областей разбиения приведенного на рис. 1, *a* и разбиения, изображенного на рис. 3, *a*, однако на 7-8 координационной окружности наблюдаются отличия в сочетаниях тайлов.

На рис. 5 приведено фрагмент квазипериодического разбиения плоскости, в котором участвует следующий набор тайлов: ромб Пенроуза; треугольник Робинсона; трапеция. Данное разбиение несколько напоминает известное разбиение Пенроуза (P3). Осуществляя гомотетию размеры новых (бо́льших) тайлов будут в τ раз больше в сравнении с размерами исходных плиток. Также из рис. 5 следует, что можно склеить по основаниям и по бедру смежные треугольники Робинсона и, следовательно, получить ромб Пенроуза (с углами 36° и 144°) и выпуклый дельтоид (с углами 72°, 72°, 72°, 144°) соответственно. Данный дельтоид является также одним из тайлом в разбиении Пенроуза (P2). Однако такое склеивание приводит к увеличению числа различных тайлов в разбиении на рис. 5 до 4. На рис. 5 можно заметить на всем фрагменте разбиения пересекающиеся декагоны.

Михалюк в [10] установил, что декагональное покрытие в области фронтов № 6, 7, 8 претерпевает топологический фазовый переход

(структурный фазовый переход Ландау 2-го рода) со спонтанным нарушением симметрии $5 \to 10$.

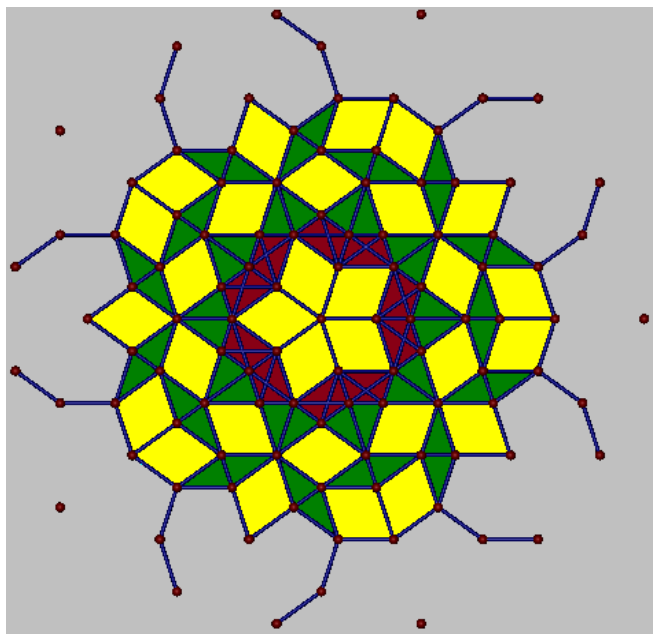

Рис. 5. Квазипериодическое разбиение плоскости (тремя различными плитками), обладающее поворотной симметрией пятого порядка.

Такой фазовый переход можно объяснить несколько в другом контексте, а именно, в вырождении пентагона в десятиугольник. В некоторых разбиениях, исследованных в данной работе, также наблюдался аналогичный топологический фазовый переход. Если принять во внимание теорию строения квазикристаллов базирующуюся на модели среза и проекции кристаллов [11], то становится понятным, чем обусловлен данный топологический фазовый переход.

Очевидно, что некоторые плоские сечения, например квазикристалла (погруженного в трехмерное пространство) ограненного, скажем в виде пентагонального додекаэдра, будут иметь симметрию как пентагона, а некоторые декагона. Все зависит от того где будет проходить срез. Если срезать вдоль определенных плоскостей симметрии, то очевидно, что симметрия может быть десятого порядка. Если же взять слой атомов, лежащий на поверхности пентагональной грани, то очевидно симметрия будет пятого порядка. Если же провести срез вдоль плоскости, лежащей между этими двумя плоскостям и параллельной им, то, как раз там мы и будем наблюдать топологический фазовый переход 5→10.

## IV. ЗАКЛЮЧЕНИЕ

В данной работе особое внимание уделено компьютерным реализациям квазипериодических разбиений плоскости, обладающих диэдральной группой симметрии $D_5$. Различных вариантов квазипериодических разбиений плоскости (с симметрией пятого порядка), в которых участвует следующий набор тайлов: треугольник Робинсона (с углами 36°, 72°, 72°), ромб Пенроуза (с углами 72° и 108°), трапеция (с углами 72° и 108°) и два подобных пентагона (коэффициент гомотетии между которыми равен $\tau$) было получено более 50 разновидностей. Отметим, что были найдены такие варианты разбиений, которые имели идентичную схожесть, вблизи и на некотором расстоянии от особой точки, но на больших расстояниях от особой точки происходило различие в сочетаниях тайлов участвующих в таких разбиениях. Другими словами, бывают такие квазипериодические разбиения, у которых одинаковый набор тайлов, которые идентичны до некоторого радиального расстояния от особой точки, но различия в сочетаниях наступают на больших расстояниях от их геометрического центра.

Методом среза и проекции были получены квазипериодические разбиения, сохраняющие всюду симметрию пятого порядка, а также и такие которые претерпевали топологический фазовый переход (структурный фазовый переход Ландау 2-го рода) со спонтанным нарушением симметрии $5 \to 10$.

Установлено, что отношение количества конкретных плиток в каждом квазипериодическом разбиении стремится к определенным значениям числа $\tau$ в некоторой степени. В разбиении Пенроуза (P3) отношение сумы числа толстых ромбов и числа тонких ромбов к числу толстых ромбов стремится к $\tau$, а отношение сумы числа толстых ромбов и числа тонких ромбов к числу тонких стремится к $\tau^2$. Если взять отношение чисел в сплаве $Al_{86}Mn_{14}$, а именно $86/14 \approx 6{,}1428571428571428$, то данный результат близок численно к $\tau^4 \approx 6{,}8541019662496845$ (заметим, что $87/13 \approx 6{,}6923076923076923$ еще ближе к $\tau^4$, однако, как известно $\tau^4$ равно отношению $n$-го члена в последовательности Фибоначчи к $(n-4)$ члену, при больших $n$). Таким образом, выдвигается гипотеза: относительное число плиток в конкретной модели квазипериодического разбиения и число их разновидностей можно сопоставить с процентным соотношением и числом компонентов в некотором сплаве соответственно, т.е. с этими модельными числами связан состав сплава, в котором возможно существование данной квазикристаллической фазы. Поскольку в каждой вершине квазипериодического разбиения плоскости можно разместить атомы (причем разного сорта), то внутренние углы тайлов, сходящихся в данных вершинах, соответствуют углам между связями для конкретных атомов с атомами из первых координационных сфер. Углы между связями характеризуют некоторые валентные углы этих атомов с учетом соответствующих гибридизаций электронных оболочек. Приведенные сведения могут быть также полезны при выборе и расчете концентраций компонентов в сплавах, в которых возможно появление квазикристаллов.

# V. СПИСОК ЛИТЕРАТУРЫ

1. Penrose R. Pentaplexity // Evreka, 1978. V. 39. p. 16-22.

2. Shechtman D. et al. // Physical Review Letters. 1984. V. 53. P. 1951.

3. Дышлис А. А., Плахтиенко Н. П. Модели нанокристаллов и неклассические периодические функции. – LAP: Lambert Academic Publishing, Deutschland, 2014. 303 с.

4. Покась С. М. Дышлис А. А, Геометрия Лобачевского и ее приложенич в математике и кристаллографии. – LAP: Lambert Academic Publishing, Deutschland, 2017. 689 с.

5. Матвеев С. В., Фоменко А. Т. Алгоритмические и компьютерные методы в трехмерной топологии. – М.: Изд-во МГУ, 1991. 301 с.

6. Фоменко А. Т. Наглядная геометрия и топология: Математические образы в реальном мирею – 2-е изд. – М.: Изд-во Моск. Ун-та, Изд-во "ЧеРо", 1998. 416 с.

7. Терстон У. Трехмерная топология и геометрия / Перевод с англ. Под ред. О.В. Шварцмана. – М.: МЦНМО, 2001. 312 с.

8. Корепин В. В. Узоры Пенроуза и квазикристаллы. // Квант. 1987. №6. с. 2–6.

9. Prokhoda A. S. About crystal lattices and quasilattices in Euclidean space. Crystallography Reports. 2017. V. 62. No. 4. p. 505–510.

10. Михалюк А. Н. Моделирование структуры декагональных квазикристаллов AL-переходной метал. 01.04.07. – физика конденсированного состояния. Владивосток, 2013. 27 с.

11. Артамонов В. А., Словохотов Ю. Л. Группы и их приложения в физике, химии, кристаллографии. – М.: Издательский центр «Академия», 2005. 512 с.